\documentclass{article}
\usepackage{amssymb, verbatim}
\usepackage{hyperref}
\usepackage{amsmath}

\def\Ric{\mathop{\rm Ric}}
\def\cRic{\mathop{\stackrel{\circ}{\Ric}}}

\def\dist{\mathop{\rm dist}}

\def\Riem{\mathop{\rm Rm}}

\def\Diam{\mathop{\rm Diam}}

\def\Vol{\mathop{\rm Vol}}
\def\VR{\mathop{\rm VR}}
\def\supp{\mathop{\rm supp}}

\def\diam{\mathop{\rm diam}}

\def\RR{\mathop{\mathbb{R}}}

\def\SS{\mathop{\mathbb{S}}}
\def\SSS{\mathop{\mathcal{S}}}

\def\be{\begin{eqnarray}}
\def\ee{\end{eqnarray}}
\def\beg{\begin{eqnarray*}}
\def\ees{\end{eqnarray*}}

\newcommand{\qed}{\hfill$\Box$}
\newtheorem{theorem}{Theorem}[section]
\newtheorem{proposition}[theorem]{Proposition}
\newtheorem{lemma}[theorem]{Lemma}

\newtheorem{corollary}[theorem]{Corollary}

\date{\small\it April 2, 2008}
\author{Brian Weber\footnote{Mathematics department, Stony Brook University, Stony Brook NY, 11794-3651. email: brweber@math.sunysb.edu}}
\title{Convergence of compact Ricci solitons}

\begin{document}

\maketitle
\begin{abstract}
We show that sequences of compact gradient Ricci solitons converge to complete orbifold gradient solitons, assuming constraints on volume, the $L^{n/2}$-norm of curvature, and the auxiliary constant $C_1$.
The strongest results are in dimension 4, where $L^2$ curvature bounds are equivalent to upper bounds on the Euler number.
We obtain necessary and sufficient conditions for limits to be compact.
\end{abstract}

\section{Introduction}\label{SectionIntroduction}

A Ricci soliton $(\SSS,g,X)$ is a Riemannian manifold $(\mathcal{S},g)$ that satisfies
\be%
2\Ric \,+\,L_Xg \,+\,2\lambda{g} &=&0, \label{EqnDefining}
\ee%
where $L$ indicates Lie differentiation, $X$ is a given vector field, and $\lambda$ is a constant.
A Ricci soliton is called gradient if $X=\nabla\,{f}$ for a function $f:\SSS\rightarrow\RR$.
It is known that compact non-Einstein gradient solitons must have $\lambda<0$; these are called shrinking gradient solitons, as opposed to the steady $\lambda=0$ and expanding $\lambda>0$ cases, which can occur on open manifolds.
It is also known that $R>0$ on shrinking solitons ($R$ is scalar curvature) \cite{Ivey1}.

Equation (\ref{EqnDefining}) and a Bianchi identity imply $2\left<\nabla^2f,\nabla{f}\right>+\nabla{R}+2\lambda\nabla{f}=0$ , so to any gradient soliton we assign the number
\be
&&C_1\;=\;|\nabla{f}|^2\,+\,R\,+\,2\lambda{f}.\label{EqnEccentricityEqn}
\ee
This is sometimes called the auxiliary equation, and $C_1$ the auxiliary constant.
Of course a normalization for $f$ is necessary before $C_1$ an have an absolute meaning, and we choose $\int_{\SSS}{f}=0$.
The constant $C_1$ has an interesting interpretation in the steady, simply connected K\"ahler case, as the divergence of a canonical holomorphic volume form \cite{Bry},
and in the context of Ricci flow, as Perelman's entropy.
Integrating (\ref{EqnEccentricityEqn}) shows $C_1\ge-n\lambda$ with equality if and only if $\SSS$ is Einstein (cf lemma \ref{LemmaW22Bounds}), so $C_1$ can perhaps be considered a measure of how far $\SSS$ is from being Einstein.

One expects the geometry of compact shrinking solitons to resemble that of positive Einstein manifolds.
Our goal is to prove a convergence-compactness result, similar to that proved for Einstein manifolds in \cite{And1}, \cite{BKN}.
In those works it was shown that the space of $n$-dimensional Einstein manifolds, with controlled volume, controlled $L^{n/2}$-norm of curvature, controlled diameter, and Einstein constants one of $\{-1,\,0,\,1\}$, is Gromov-Hausdorff precompact with limiting objects being Einstein manifolds with point-singularities of $C^\infty$ orbifold type.

In the case of Einstein constant $1$, Myer's theorem implies diameter bounds from the other hypotheses, and in 4-dimensions the $L^2$-norm of curvature is just the Euler characteristic.
Thus the 4-dimensional positive Einstein case involves attractively few assumptions: a sequence of such manifolds with lower volume bounds and upper bounds on the Euler characteristic degenerates at worst to a compact Einstein orbifold with (uniformly) finitely many $C^\infty$-type orbifold points.

In the soliton case, one expects $C_1$ to play a role as well.
Noncompact shrinking solitons exist (e.g. the Gaussian soliton), so volume upper bounds should be expected.
Let $\mathfrak{S}^n=\mathfrak{S}^n(\underline{V},\overline{V},\overline{C_1},\Lambda)$ be the set of compact, $n$-dimensional gradient Ricci solitons $(\SSS^n,g,\nabla{f})$ that satisfy
\begin{itemize}
\item[{\it i})] the normalizations $\lambda=-1$ and $\int_{\SSS}f=0$,
\item[{\it ii})] upper and lower volume bounds $0<\underline{V}<\Vol\SSS<\overline{V}<\infty$,
\item[{\it iii})] upper bounds $C_1<\overline{C_1}<\infty$ on $C_1$, and
\item[{\it iv})] upper bounds on $L^{\frac{n}{2}}(|\Riem|)$: $\int_{\SSS}|\Riem|^{\frac{n}{2}}<\Lambda$.
\end{itemize}
Diameter bounds do not seem to follow from the other assumptions, but there is a natural choice of basepoints for taking pointed limits.
Given a soliton $(\SSS,g,\nabla{f})\in\mathfrak{S}$, let $p_m\in\SSS$ be a point with $f(p_m)=\min_{\SSS}f$.
The distance between any two such points is bounded in terms of $C_1$ and $\lambda$ (lemma \ref{LemmaQuadraticf}), so it does not matter which point in $f^{-1}(\min{f})$ is chosen.
This choice is also natural in that the unit ball centered at $p_m$ is noncollapsed (theorem \ref{ThmOneNoncollapsedBall}).

\begin{theorem}\label{ThmMainNDimTheorem}
Let $(\SSS_i^n,g_i,\nabla{f}_i)$ be a sequence of solitons in $\mathfrak{S}^n$, each with basepoint $p_i$.
Then a subsequence converges in the pointed Gromov-Hausdorff topology to a complete (possibly noncompact) orbifold gradient shrinking soliton $(\SSS_\infty,p_\infty)$ with locally finitely many singular points.
The singular points are all orbifold points, and away from the singularities the convergence is locally in the $C^\infty$ topology.
\end{theorem}
In the odd-dimensional case, it is impossible for singularities to form a finite distance from the basepoint, due to the fact that singularity models are oriented, Ricci-flat, non-flat, and ALE and so do not exist in odd dimensions (this is pointed out, for instance, in the proof of theorems $A$ and $A'$ in \cite{And1}).

On 4-dimensional compact solitons, the $L^2$-norm of curvature is controlled by the Euler number and the (scale invariant) number ${C_1}^2\Vol\SSS$.
Let $\mathfrak{S}^4_\chi=\mathfrak{S}^4_\chi(\underline\lambda,\overline\lambda,\overline{C_1},\chi)$ be the set of compact, 4-dimensional gradient Ricci solitons $(\SSS,g,\nabla{f})$ with condition ({\it iv}) replaced by
\begin{itemize}
\item[{\it iv'})] Upper bounds on the Euler characteristic: $\chi(\SSS)<\chi$.
\end{itemize}
Lower bounds are unnecessary because $\chi(\SSS)\ge2$, which follows from results in \cite{Lott} \cite{Derd1} \cite{FLGR} where the finiteness of the fundamental group on compact shrinking solitons is proved (something (probably) slightly weaker is proved in \cite{Derd1}).
See also \cite{FMZ} for the noncompact, bounded Ricci-curvature case.
\begin{theorem}\label{ThmBoundingRmIn4D}
On a 4-dimensional compact gradient Ricci soliton,
\beg
\int|\Riem|^2&=&8\pi^2\chi\,+\,2\lambda^2\Vol\SSS\,+\,\frac38\int\left(R-\overline{R}\right)^2\\
&\le&8\pi^2\chi\,+\,\frac{27}{128}C_1^2\Vol\SSS.
\ees
where $\overline{R}=\frac{1}{\Vol\SSS}\int_{\SSS}R=-4\lambda$.
\end{theorem}
\begin{corollary} \label{Cor4DimLimits}
Theorem \ref{ThmMainNDimTheorem} holds for the class $\mathfrak{S}^4_\chi$.
Additionally, any limiting object $M_\infty$ has only uniformly finitely many singularities, and if it is noncompact, it is collapsed with locally bounded curvature at infinity.
\end{corollary}
The phrase ``collapsed with locally bounded curvature at infinity'' means that given any $\kappa>0$, there is a compact set $K\subset\subset{M}_\infty$ that $M_\infty-K$ is $\kappa$-collapsed at every point, where, following Perelman, $\kappa$-collapsed at $p\in{M}$ means that if $\sup_{B(p,r)}|\Riem|=r^{-2}$, then $\Vol{B}(p,r)\le\kappa{r}^n$. This is essentially the same as other definitions in the literature, eg \cite{And4}.

We obtain a description of when limits of theorem \ref{ThmMainNDimTheorem} are compact.
\begin{theorem}\label{ThmUpperBoundsImplyUpperBounds}
Given an n-dimensional compact shrinking soliton $\SSS$ with $\lambda=1$ and definite bounds on $\Vol\SSS$ and $C_1$, then the value of any one of the following quantities implies definite bounds on any of the others:
\beg
&&\Diam\,{\SSS} \quad\quad\quad
{\sup}_{\SSS}{R} \quad\quad\quad
{\sup}_{\SSS}{|\nabla{f}|} \quad\quad \quad
{\sup}_{\SSS}{f} \quad\quad\quad
{\inf}_{p\in\SSS}\Vol\,B_p(1).
\ees
Additionally, upper bounds on the Sobolev constant imply any of the above.
In the 4-dimensional case, $\kappa$-noncollapsing is also equivalent to any of the above.
\end{theorem}

Moduli spaces of Ricci solitons have been studied in \cite{CS} and \cite{Zhang}, upon whose results we build.
Our improvements are the relaxation of their diameter and pointwise Ricci curvature bound hypotheses.


{\bf Remark}. It is known that every compact Ricci soliton is gradient, for possibly a different vector field $X$. This was proved by Perelman \cite{Perel1}; see also \cite{Cao2} for a concise proof based on Perelman's work.
This is true for open, shrinking solitons with bounded curvature \cite{Nab}, but not true for expanding solitons (eg \cite{Lau}).
It has been known for some time that compact 2- and 3-dimensional Ricci solitons are Einstein \cite{Ham1} \cite{Ivey1}.

{\bf Remark}. This work partially answers in the affirmative a conjecture of Cheeger-Tian's \cite{CT}, but we do not use their methods, nor achieve the full strength of their results.
Due to collapsing phenomena as diameters grow, we cannot rule out that unboundedly many singularities develop, even in the 4-dimensional case, although the only possible accumulation point is at infinity.

{\bf Organization}. In section \ref{SectionRicciSolitons} we collect information on the functions $f$ and $|\nabla{f}|$, establish $L^2(|\Riem|)$ bounds in the 4-dimensional case, and prove most of theorem \ref{ThmUpperBoundsImplyUpperBounds}.
In section \ref{SectionVolumeComparison} we write down a version of relative volume comparison for gradient solitons, and prove the existence of a noncollapsed region in any compact shrinking soliton as well as volume ratio control on all balls finite distances away from the basepoint.
We also complete the proof of theorem \ref{ThmUpperBoundsImplyUpperBounds}.
In section \ref{SectionEllipticRegularity} we run through the $\epsilon$-regularity theorem, and establish local bounds on Sobolev constants under certain conditions.
We prove theorem \ref{ThmMainNDimTheorem} in section \ref{SectionConvergence}.

{\bf Acknowledgements}. I would like to thank Prof. Xiuxiong Chen for suggesting this problem, and for his encouragement.
I would also like to thank Bing Wang for some very useful suggestions.

\section{Compact Ricci Solitons}\label{SectionRicciSolitons}

Here we collect some information on Ricci solitons which can be verified without $\epsilon$-regularity.
With one exception these results apply in any dimension.
Throughout this section we normalize so $\int{f}=0$.
First we identify the elliptic equations necessary $\epsilon$-regularity.
These equations are very standard, so our proof is brief.
\begin{proposition}[Elliptic Equations]\label{PropEllipticResults}
If $(\SSS,g,\nabla{f})$ is any $n$-dimensional gradient Ricci soliton, then
\be
&&\triangle{f}\;=\;-R-n\lambda \label{EqnTrianglef}\\
&&\triangle{d}f\;=\;-\Ric(\nabla{f}) \label{EqnTriangleNablaf}\\
&&dR\;=\;2\delta\Ric\;=\;2\Ric(\nabla{f}) \label{EqnDivRic}\\
&&\delta\Riem\left(\cdot,\cdot,\cdot\right)\;=\;\Riem\left(\nabla{f},\cdot,\cdot,\cdot\right) \label{EqnDivRiem}\\
&&\triangle\Riem\;=\;\nabla(\Riem*\nabla{f})\,+\,\Riem*\Riem \label{EqnTriangleRiem}\\ &&\triangle\Ric\;=\;\nabla^2{R}\,+\,\Riem*\Ric\,+\,\lambda\Ric\,+\,\Riem*(\nabla{f}\otimes\nabla{f}) \label{EqnTriangleRic}\\
&&\triangle{R}\;=\;-2|\Ric|^2\,-\,2\lambda{R}\,+\,\left<\nabla{R},\nabla{f}\right> \label{EqnTriangleR}
\ee
\end{proposition}
\underline{\sl Pf}\\
\indent Equation (\ref{EqnTrianglef}) is just the traced soliton equation.
We have
\beg
{\Ric}_{ij,j}&=&-f_{,ijj}\;=\;-f_{,jji}-{\Ric}_{is}f_s\;=\;R_{,i}-{\Ric}_{is}f_s.
\ees
Using $\Ric_{ij,j}=\frac12R_{,i}$ we get (\ref{EqnDivRic}).
Equation (\ref{EqnTriangleNablaf}) also follows easily.
Then (\ref{EqnDivRiem}) is a Bianchi identity:
\beg
{\Riem}_{sijk,s}&=&{\Ric}_{ij,k}-{\Ric}_{ik,j}\;=\;f_{,ikj}-f_{,ijk}\;=\;{\Riem}_{sijk}f_{,s}.
\ees
We omit the straightforward computations for (\ref{EqnTriangleRiem}) and (\ref{EqnTriangleRic}).
The exact expression for (\ref{EqnTriangleRic}) is
\beg
{\Ric}_{ij,ss}&=&\frac12R_{,ij}\,+\,2\,{\Riem}_{sijk}{f}_{sk}\,-\,{\Ric}_{js}{f}_{is}\,+\,{\Riem}_{sijk}f_{,s}f_{,k}
\ees
so that (\ref{EqnTriangleR}) follows by tracing.
\qed

Let $\overline{R}=\frac{1}{\Vol\SSS}\int_{\SSS}R=-n\lambda$ denote the average value of $R$. If $\triangle$ is the ``negative'' Laplacian, we can write for instance $\triangle{f}=-\left(R-\overline{R}\right)$.
Notice an immediate conclusion of the following lemma is that $C_1>n\lambda$ if the soliton is not Einstein.
\begin{lemma}[$W^{1,2}$-bounds on $|\nabla{f}|$]\label{LemmaW22Bounds}
On any compact gradient Ricci soliton $\SSS$,
\be
&&\int_{\SSS}|\nabla^2f|^2\;=\;\frac12\int_{\SSS}\left(R\,-\,\overline{R}\right)^2\label{EqnL2Nabla2f}\\
&&\int_{\SSS}|\nabla{f}|^2\;=\;\left(C_1\,+\,n\lambda\right)\Vol\SSS\label{IneqL2Nablaf}
\ee
\end{lemma}
\underline{\sl Pf}\\
\indent For (\ref{EqnL2Nabla2f}),
\beg
\int|\nabla^2f|^2&=&\int\left(\triangle{f}\right)^2\,-\,\int\Ric(\nabla{f},\nabla{f})\;=\;\int\left(\triangle{f}\right)^2\,-\,\frac12\int\left<\nabla{R},\nabla{f}\right>\\
&=&\int\left(\triangle{f}\right)^2\,+\,\frac12\int{R}\triangle{f}\;=\;\frac12\int\left(\triangle{f}\right)^2.
\ees
For (\ref{IneqL2Nablaf}) we integrate the ``auxiliary equation'' to get
\beg
\int|\nabla{f}|^2&=&\int\left(-R\,+\,C_1\right)\;=\;\int(n\lambda\,+\,C_1).
\ees
\qed

\begin{proposition}[$L^2$ bounds on curvature]\label{PropL2CurvatureBounds}
On a compact gradient Ricci soliton,
\beg
&&\int_{\SSS}R^2\;\le\;-\lambda\left((n+2)C_1\,+\,2n\lambda\right)\Vol\SSS\\
&&\int_{\SSS}(R\,-\,\overline{R})^2\;\le\;-\lambda\left(C_1+n\lambda\right)(n+2)\Vol\SSS\\
&&\int_{\SSS}|\Ric|^2\;=\;-\frac12\lambda^2n(n-2)\Vol\SSS\,+\,\frac12\int{R}^2\\
&&\int_{\SSS}|\cRic|^2\;=\;\frac{n-2}{2n}\int_{\SSS}\left(R\,-\,\overline{R}\right)^2.
\ees
On a 4-dimensional compact gradient Ricci soliton,
\beg
\int|\Riem|^2&=&8\pi^2\chi\,+\,2\lambda^2\Vol\SSS\,+\,\frac38\int\left(R-\overline{R}\right)^2
\ees
\end{proposition}
\underline{\sl Pf}\\
\indent We first find an expression for $\int{R}|\nabla{f}|^2$.
\beg
\int{R}|\nabla{f}|^2&=&-\int\triangle{f}|\nabla{f}|^2\,-\,n\lambda\int|\nabla{f}|^2\\
&=&2\int\nabla^2{f}(\nabla{f},\nabla{f})\,-\,n\lambda\int|\nabla{f}|^2\\
&=&-2\int|\nabla^2f|^2\,-\,\lambda(n+2)\int|\nabla{f}|^2.
\ees
Since $R|\nabla{f}|^2\ge0$ it follows that
\beg
\int\left(R-\overline{R}\right)^2&\le&-\lambda(n+2)\int|\nabla{f}|^2\\
\int{R}^2&\le&\overline{R}^2\Vol\SSS\,-\,\lambda(n+2)\left(C_1+n\lambda\right)\Vol\SSS\\
&=&-\lambda\left((n+2)C_1\,+\,2n\lambda\right)\Vol\SSS
\ees
The rest now follows from
\beg
\int{R}^2&=&n^2\lambda^2\Vol\SSS\,+\,\int\left<\nabla{f},\nabla{R}\right>\;=\;n^2\lambda^2\,+\,2\int\left<\delta\Ric,\nabla{f}\right>\\
&=&n^2\lambda^2\Vol\SSS\,+\,2\int|\Ric|^2\,+\,2\lambda\int{R}
\ees
and the usual formulas
\beg
8\pi^2\chi(\SSS)&=&\frac{1}{24}\int{R}^2\,-\,\frac12\int|\cRic|^2\,+\,\int|W|^2\\
\int_{\SSS}|\Riem|^2&=&\frac16\int{R}^2\,+\,\frac12\int|\cRic|^2\,+\,\int|W|^2.
\ees
\qed

\begin{proposition}[Growth upper bounds on $f$]\label{PropUpperFBounds}
If $(\SSS^n,g,\nabla{f})$ is any shrinking Ricci soliton, then given $p,q\in\SSS$ we have
\beg
f(q)&\le&2f(p)\,+\,\frac{C_1}{2|\lambda|}\,+\,|\lambda|\left(\dist(p,q)\right)^2
\ees
\end{proposition}
\underline{\sl Pf}\\
\indent Since $R>0$ we have
\beg
&&|\nabla{f}|^2\;=\;C_1\,-\,R\,-\,2\lambda{f}\;<\;C_1\,-\,2\lambda{f}\\
&&|\nabla(C_1-2\lambda{f})^{1/2}|\;<\;|\lambda|.
\ees
Along a unit-speed geodesic between $p$ and $q$, this differential inequality has the solution
\beg
&&\sqrt{C_1-2\lambda{f}(p)}\;\le\;\sqrt{C_1-2\lambda{f}(q)}\,+\,|\lambda|\dist(p,q).
\ees
Squaring both sides and using a Schwartz inequality gives the result.
\qed

Proposition \ref{PropUpperFBounds}, showing that $f$ grows at most quadratically, holds whether the soliton is compact or not. Lemma \ref{LemmaQuadraticf} below shows that the growth of $f$ is bounded quadratically from below as well.
Comparing this to lemma 2 of \cite{Bry} for the steady case and theorem 1.5 of \cite{Cao1} for the expanding case shows that this growth rate is quite characteristic of shrinking solitons.
Choosing $q$ to be any point with $f(q)\le0$, this inequality now bounds $f$, and therefore $R$, in terms of diameter. This is the first half of the following proposition.

\begin{proposition}\label{PropDiamRBonds}
Given constraints on $C_1$ and $\lambda$ on a compact gradient Ricci soliton, an upper bound on diameter is equivalent to an upper bound on scalar curvature. Specifically,
\be
{\max}_{\SSS}{R}&\le&2C_1\,+\,2\lambda^2(\diam\,\SSS)^2,\label{IneqRByDiam}\\
{\diam}\,\SSS&\le&\frac{\pi}{|\lambda|}\sqrt{2(n-1)|\lambda|+C_1\,+\,{\max}_{\SSS}{R}}\label{IneqDiamByR}.
\ee
\end{proposition}
\underline{\sl Pf}\\
\indent Choosing $p$ so that $f(p)<0$, inequality (\ref{IneqRByDiam}) is immediate from proposition \ref{PropUpperFBounds}.

If $R$ is bounded, then $\max{f}<(\max{R})/2|\lambda|$ implies $|\nabla{f}|^2<C_1+\max{R}$.
This allows us to follow the standard Myer's argument for bounding diameters. Letting $\gamma(r)$ be any unit-speed, minimizing geodesic of length $d$ and $g(r)$ any 1-variable function with $g(0)=g(d)=0$, a variational argument provides the standard inequality
\beg
0&\le&\int_0^d(n-1)\left(g'(r)\right)^2\,-\,\left(g(r)\right)^2\Ric(\dot\gamma,\dot\gamma).
\ees
Substituting the expression for $\Ric$ and using integration by parts,
\be
0&\le&\int_0^d(n-1)\left(g'(r)\right)^2\,+\,\lambda\left(g(r)\right)^2\,-\,2g(r)g'(r)f'(r)\label{IneqVariationOfEnergy}\\
0&\le&\int_0^d\left(n-1-\frac2\lambda|\nabla{f}|^2\right)\left(g'(r)\right)^2\,+\,\frac{\lambda}{2}\left(g(r)\right)^2\nonumber
\ee
(recalling that $\lambda$ is negative). It is standard to choose $g(r)=\sin(\pi{r}/d)$; we get
\beg
0&\le&\frac{\pi}{2d}\left(n-1-\frac2\lambda\max{|\nabla{f}|^2}\right)\,+\,\frac{\lambda}{2}\frac{d}{2\pi}.
\ees
This bounds the length $d$ in terms of $\lambda$, $C_1$, and $\max{R}$.
\qed

In fact, this variational argument provides a lower growth estimate for $f$.
\begin{lemma}[Growth lower bounds on $f$]\label{LemmaQuadraticf}
On any shrinking soliton it holds that
\be
f(p)\,+\,f(q)\;\ge\;\frac{-\lambda}{16\pi}\dist(p,q)^2\,-\,(n-1)\pi\,+\,C_1/\lambda
\ee
\end{lemma}
\underline{\sl Pf}\\
\indent Let $\gamma$ be a unit-speed, minimizing geodesic of length $l$, parameterized by the variable $r$. We use (\ref{IneqVariationOfEnergy}) in the form
\beg
0&\le&\int_0^l(n-1)\left(g'(r)\right)^2\,+\,\lambda\int_0^l\left(g(r)\right)^2 \,+\,\int_0^l(g(r)^2)\,f''(r)
\ees
Picking $s>0$ with $2s<l$, let $g(r)$ be the $C^{1,1}$ function
\beg
g(r)&=&\begin{cases}
\sin(\pi\,r/2s)     & 0\le{r}<s\\
1                   & s\le{r}<{l-s}\\
\sin(\pi\,(l-r)/2s) & {l-s}\le{r}\le{l}.
\end{cases}
\ees
Using integration by parts twice on the last term, we get
\beg
0&\le&\frac{\pi^2(n-1)}{2s}\,+\,\lambda\,\left(l\,-\,s\right) \,+\,\frac{\pi^2}{2s^2}\int_0^s\cos\left(\frac{\pi{r}}{s}\right)\left(f(r)+f(l-r)\right)\,dr
\ees
But by proposition \ref{PropUpperFBounds}, we have the bounds
\beg
-\frac{C_1}{|\lambda|}\;\le\;f(r)+f(l-r)\;\le\;2f(0)+2f(l)\,+\,\frac{C_1}{|\lambda|}\,+\,2|\lambda|\,r^2,
\ees
so we can estimate the last integral in terms of $f(0)+f(l)$ only. We get
\beg
0&\le&\frac{\pi}{s}\left(\frac{\pi(n-1)}{2}\,+\,f(0)+f(l)+\frac{C_1}{|\lambda|}\right)\,+\,\lambda\,\left(l\,-\,4s\right).
\ees
Putting $s=l/8$ say, we conclude
$\frac{|\lambda|}{16\pi}\,l^2 \,-\,\frac{\pi(n-1)}{2}\,-\,\frac{C_1}{|\lambda|}\,<\,f(0)+f(l).$
\qed

\section{Volume comparison for Ricci solitons}\label{SectionVolumeComparison}

\begin{lemma}\label{LemmaLaplacianComparison}
Let $r=\dist(p,\cdot)$ denote the distance to some point $p\in\SSS$. There exists a function $F:\RR^{+}\rightarrow\RR^{+}$, which depends only on $n$, $C_1$, $\lambda$, and $f(p)$, so that
\beg
\triangle{r}&\le&F(r).
\ees
We can choose $F$ so that $\lim_{r\rightarrow0^{+}}r\cdot{F}(r)\;=\;n-1$.
Note that $r^{-n}\cdot{F}(r)\rightarrow\omega_n$ as $r\searrow0$.
\end{lemma}
\underline{\sl Pf}\\
\indent If $r$ is the distance function to the point $p$, then as usual
\beg
\frac{d}{dr}\triangle{r}\,+\,|\nabla^2r|^2\,+\,\Ric\left(\frac{d}{dr},\frac{d}{dr}\right)&=&0,
\ees
which, combined with the soliton equation gives
\beg
&&\frac{d}{dr}\left(\triangle{r}\,-\,\left<X,\nabla{r}\right>\right)\,+\,\frac{\left(\triangle{r}\right)^2}{n-1}\,-\,\lambda\;\le\;0.
\ees
Note that in any case
\beg
&&\frac{d}{dr}\left(\triangle{r}\,-\,\left<X,\nabla{r}\right>\right)\;\le\;\lambda.
\ees
If $\frac{\alpha}{1-\alpha}(\triangle{r})^2\ge(\left<X,\nabla{r}\right>)^2$ for some $0<\alpha<1$, then
\be
&&\frac{d}{dr}\left(\triangle{r}\,-\,\left<X,\nabla{r}\right>\right)\,+\,\frac{(1-\alpha)(\triangle{r})^2+\alpha(\triangle{r})^2}{n-1}\,-\,\lambda\;\le\;0\nonumber\\
&&\frac{d}{dr}\left(\triangle{r}\,-\,\left<X,\nabla{r}\right>\right)\,+\,\frac{1-\alpha}{2}\frac{\left(\triangle{r}-\left<X,\nabla{r}\right>\right)^2}{n-1}\,-\,\lambda\;\le\;0.\label{IneqSolitonBishop}
\ee
Integrating this differential inequality with appropriate (asymptotic) initial conditions gives
\be
\triangle{r}&\le&\left<X,\nabla{r}\right>\,+\,\sqrt{\frac{-2(n-1)\lambda}{1-\alpha}}\,\cdot\,\cot\left(r\cdot\sqrt{\frac{-\lambda(1-\alpha)}{2(n-1)}}\right),\label{IneqCompIneq1}
\ee
which holds until $\frac{\alpha}{1-\alpha}(\triangle{r})^2=(\left<X,\nabla{r}\right>)^2$.

Assume this holds true at least until distance $r=r_0$ away from $p$. Note that in the gradient case $X=\nabla{f}$, because of the growth condition on $f$, we can bound $r_0$ uniformly from below, depending on the values of $\lambda$, $C_1$, and $\lambda{f}(p)$. When $r>r_0$ we can still say that
\beg
\triangle{r}&\le&\left<X,\nabla{r}\right>\,+\,\triangle{r}\big|_{r=r_0}\,-\,\left<X,\nabla{r}\right>\big|_{r=r_0}\,+\,\lambda\cdot(r-r_0)\nonumber\\
\triangle{r}&\le&\left<X,\nabla{r}\right>\,+\,\left(\frac{\sqrt{1-\alpha}}{\sqrt\alpha}\,+\,1\right)\left<X,\nabla{r}\right>\big|_{r=r_0}\,+\,\lambda\cdot(r-r_0).
\ees
If the soliton is gradient then $|X|=|\nabla{f}|\le\sqrt{C_1-2\lambda{f}(p)}-\lambda\cdot{r}$, so
\be
\triangle{r}&\le&\left(\frac{\sqrt{1-\alpha}}{\sqrt{\alpha}}\,+\,1\right)\left(\sqrt{C_1-2\lambda{f}(p)}\,+\,\lambda\cdot{r}\right).\label{IneqCompIneq2}
\ee

We get our result: with any $0<\alpha<1$ we have
\beg
\triangle{r}\;\le\;\begin{cases}
\sqrt{C_1-2\lambda{f}(p)}-\lambda\cdot{r}\,+\,\sqrt{\frac{-2(n-1)\lambda}{1-\alpha}}\,\cdot\,\cot\left(r\cdot\sqrt{\frac{-\lambda(1-\alpha)}{2(n-1)}}\right) & {\rm if}\;0<r<r_0\\
\left(\sqrt{\frac{1-\alpha}{\alpha}}\,+\,1\right)\left(\sqrt{C_1-2\lambda{f}(p)}\,+\,\lambda\cdot{r}\right) & {\rm if} \;r\ge{r}_0
\end{cases}
\ees
To see that $F$ can be chosen so that $r\cdot{F}(r)\rightarrow(n-1)$ as $r\searrow0$, just make $\alpha$ closer and to $0$, and note the asymptotic behavior of the cotangent function.
\qed

\begin{theorem}[Soliton relative volume comparison]\label{ThmSolitonVolumeComparison}
Assume $(\SSS,g,f)$ is a gradient Ricci soliton, and let $p\in\SSS$. Let $\overline{F}$ be the function $\overline{F}(r)=\omega_n\int_0^r\int_0^uF(s)s^{n-1}\,ds\,du$, where $\omega_n$ is the volume of the Euclidean unit ball. Then $0\le{r}<R$, $0\le{s}<S$, $s\le{r}$, $S\le{R}$ implies
\beg
\frac{\Vol{B}(p,R)-\Vol{B}(p,r)}{\Vol{B}(p,S)-\Vol{B}(p,s)}&\le&\frac{\overline{F}(R)-\overline{F}(r)}{\overline{F}(S)-\overline{F}(s)}.
\ees
The function $\overline{F}$ depends only on $n$, $\lambda$, $C_1$, and $f(p)$.
\end{theorem}
\underline{\sl Pf}\\
\indent The conclusion is immediate from lemma \ref{LemmaLaplacianComparison} using a standard arguments, but we briefly run through the steps involved. For any $t>0$ integration by parts gives
\beg
|\partial{B}(p,t)|&=&\int_{B(p,r)}\triangle{r}\,d\Vol\\
&=&\int_{\SS^{n-1}}\int_0^t\triangle{r}\,r^{n-1}\,dr\,d\sigma\\
&\le&\omega_n\int_0^tF(r)r^{n-1}dr.
\ees
Of course then $\Vol{B}(p,r)=\int_0^r|\partial{B}(p,s)|ds\le\overline{F}(r)$; this is commonly referred to as Bishop volume comparison. To get the full Bishop-Gromov result, one uses the fact that, for any functions $h$, $j$ with $h'\le{j}'$ and the same initial conditions, one has that the function
\beg
K(x,y)\;=\;\frac{\int_x^yh(t)\,dt}{\int_x^yj(t)\,dt}
\ees
is (nonstrictly) decreasing in both variables.
\qed

We can use this version of volume comparison along with the growth upper bounds on $f$ to establish a noncollapsing result.
\begin{theorem}\label{ThmOneNoncollapsedBall}
If $(\SSS^n,g,f)$ is a compact Ricci soliton, there is a noncollapsed unit ball. Specifically, there is a constant $v=v(n,C_1,\lambda)$ so that $\Vol{B}(p_m,1)>v$ where $p_m$ is any point realizing the minimum of $f$.
\end{theorem}
\underline{\sl Pf}\\
\indent Let $f_{+}$ and $f_{-}$ be the positive and negative parts of $f$. Since $\min{f}\ge{C}_1/2\lambda$ we have
\be
-\frac{C_1}{2\lambda}\Vol\SSS&\ge&\int{f}_{-}\;=\;\int{f}_{+}\;\ge\;\alpha\Vol\left(\{f\ge\alpha\}\right)\nonumber\\
\Vol\left(\{f\ge\alpha\}\right)&\le&\frac{C_1}{2|\lambda|\alpha}\Vol\SSS.\label{IneqSuperLevelSetControl}
\ee
By lemma \ref{LemmaQuadraticf}, given any $\alpha$ we have that $f\ge\alpha$ outside $B(p_m,\sqrt{-C(n,C_1,\lambda)+\alpha/51})$. Thus for a given large value of $r$, $B(p_m,r)$ contains most of the soliton's volume.
But soliton relative volume comparison holds inside $B(p_m,r)$, so that $\Vol{B}(p_m,1)$ is controlled from below. \qed

\begin{corollary} \label{CorNoncollapsedBalls}
Given $p\in\SSS$, there is a $V>0$ depending on $n$, $f(p)$ (or equivalently $\dist(p,p_m)$), $C_1$, and $\lambda$ so that $Vol\,B(p,r)\ge{V}\,r^n$ for $0<r\le1$.
\end{corollary}
\underline{\sl Pf}\\
\indent Because the noncollapsed ball $p_m$ is finitely far from $p$, soliton relative volume comparison for balls centered at $p$ will give the result.
\qed

\section{Elliptic Regularity}\label{SectionEllipticRegularity}

We prove two main theorems here.
The first is an improvement on the standard $\epsilon$-regularity theorems, where bounds on local Sobolev constants are replaced by bounds on local volume ratios.
This result holds for any metric satisfying the standard $\epsilon$-regularity, not just Ricci solitons.
The proof is based on proposition 3.1 of \cite{TV3}.
Secondly we prove $\epsilon$-regularity for Ricci solitons, following \cite{CS} \cite{Zhang}.

Given a Riemannian metric on an open set $\Omega$, it's Sobolev constant $C_S=C_S(\Omega)$ is the smallest number such that
\beg
\left(\int\phi^{\frac{2n}{n-2}}\right)^{\frac{n-2}{n}}\;\le\;C_S\int|\nabla\phi|^2
\ees
holds for any function $\phi\in{C}^1_c(\Omega)$. If $\diam(\Omega)<\infty$ say, then necessarily $C_S>0$. The Sobolev constant has a well-known geometric interpretation: an upper bound on $C_S(\Omega)$ is equivalent to an upper bound on the isoperimetric inequality on subdomains of $\Omega$.

Given an elliptic equation for some function $u$ and $L^{n/2}$-bounds on $u$, one can use an iteration process to establish $L^p$-bounds, or if the original bounds are $L^p$ for $p>n/2$, $L^\infty$-bounds.
The specific estimates rely on the Sobolev constant, so are useless if $C_S$ is not controlled.

The Riemann tensor of a canonical Riemannian metric (e.g. Einstein, CSC-Bach flat, extremal K\"ahler, etc.) normally satisfies some elliptic equation; in the Ricci soliton case we have proposition \ref{PropEllipticResults}.
Thus higher regularity estimates are possible given $L^p$ bounds on $|\Riem|$.
Authors normally have little use for determining exactly how the Sobolev constant appears in the estimates, but if careful track is kept, $\epsilon$-regularity would usually read
\be
C_S^{\frac{n}{2}}\int_{B(p,r)}|\Riem|^{\frac{n}{2}}\,\le\,\epsilon
\quad\Rightarrow\quad \sup_{B(p,r/2)}|\Riem|\,\le\,C\cdot{r}^{-2}\,\left(C_S^{\frac{n}{2}}\int_{B(p,r)}|\Riem|^{\frac{n}{2}}\right)^\frac2n\label{IneqStandardReg}
\ee
where $C_S$ is the local Sobolev constant in the ball $B(p,r)$, and the constants $C$ and $\epsilon$ do not depend on $C_S$.
Sobolev constants (or isoperimetric constants) can be difficult to estimate, but with an elliptic system for curvature the Sobolev constant can be made to depend on volume ratios in small-energy regions.
We determine the explicit dependency in two steps.
We use the notation
\beg
\VR\,B(p,r)\;=\;r^{-n}\Vol\,B(p,r).
\ees

\begin{lemma} \label{LemmaAlmostFlatSobolevBounds}
There are dimensional constants $\epsilon_1$, $C$ so that $\sup_{B(p,r)}|\Riem|<\epsilon_1$ implies $C_S(B(p,r))<C\cdot\left(\VR\,B(p,r)\right)^{-\frac2n}$.
\end{lemma}
\underline{\sl Pf}\\
\indent Everything is scale-invariant, so we can work in just the unit ball $B=B(p,1)$. To deal with the possibility of geometric collapsing we first pass to a noncollapsing domain as follows.
If $\epsilon_1$ is small enough, the exponential map at $p$ is noncritical on a ball of radius 1 in the tangent space.
We work in this region of the tangent space, with the pullback metric.
There is some maximal, open star-shaped domain $0\in\overline{Q}_{0}\subset{T}_pM$ on which $\exp_p$ is a diffeomorphism ($\overline{Q}_0$ is uniquely determined by requiring its image be the ball $B(p,1/2)$ minus the cut locus of $\exp_p$).
Given any point $\bar{p}\in\exp^{-1}_p(p)$ there exists a ``translate'' $\overline{Q}_{\bar{p}}$ of $\overline{Q}_0$, obtained by lifting the image of $\overline{Q}_0$ starting at the basepoint $\bar{p}\in{T}_pM$ rather than at $0\in{T}_pM$.
Now let $\overline{Q}$ be the closure of the union of all the $\overline{Q}_{\bar{p}}$ so that $\dist(\bar{p},0)<\frac12$.
Note that $\overline{Q}\subset{B}(1)\subset{T}_pM$.
We use the pullback metric on $\overline{Q}$.

Off a set of measure $0$, $\exp_p:\overline{Q}\rightarrow{B}(p,1/2)$ is precisely a $k$-$1$ map.
We estimate $k$ using Bishop volume comparison: for $|\Riem|$ small enough, then
\beg
k&=&\frac{\Vol\overline{Q}}{\Vol{B}(p,1/2)}\;\le\;\frac{2^{n+1}\Vol{B}_{Eucl.}(1)}{2^n\Vol{B}(p,1/2)}\;=\;\frac{C(n)}{\VR\,B(p,1/2)}.
\ees

Our final requirement on $\epsilon_1$ is that it be small enough that on any ball of radius $1$ which is diffeomorphic to a Euclidean ball and has $|\Riem|<\epsilon_1$, the Sobolev constant is less than twice the Euclidean Sobolev constant, $C_E$.

Now let $f\in{C}^1(B)$, and lift it to a function $\bar{f}$ on $\overline{Q}$ (it will still be $C^1$). Since $|\Riem|<\epsilon_1$ on $\overline{Q}$, an equivalent version of the Sobolev inequality gives
\beg%
\left(\int_{\overline{Q}}|\bar{f}|^{\frac{2n}{n-2}}\right)^{\frac{n-2}{n}} &\le&2C_E\int_{\overline{Q}}|\nabla\bar{f}|^2 \,+\,\left(\Vol\overline{Q}\right)^{-\frac{2}{n}}\int_{\overline{Q}}|\bar{f}|^2.
\ees%
Therefore
\beg%
\left(k\int_{B}|\bar{f}|^{\frac{2n}{n-2}}\right)^{\frac{n-2}{n}} &\le&2C_Ek\int_{B}|\nabla{f}|^2\,+\,\left(k\Vol{B}\right)^{-\frac{2}{n}}k\int_{B}|f|^2\\
\left(\int_{B}|f|^{\frac{2n}{n-2}}\right)^{\frac{n-2}{n}} &\le&2C_Ek^{\frac2n}\int_{B}|\nabla{f}|^2 \,+\,\left(\Vol{B}\right)^{-\frac{2}{n}}\int_{B}|f|^2.
\ees%
Thus finally
\beg%
\left(\int_{B}|f|^{\frac{2n}{n-2}}\right)^{\frac{n-2}{n}} &\le&C_1\cdot{VR}(B)^{-\frac2n}\int_{B}|\nabla{f}|^2 \,+\,\left(\Vol{B}\right)^{-\frac{2}{n}}\int_{B}|f|^2.
\ees%
for a dimensional constant $C_1$. \qed

\begin{proposition}\label{PropVolEpsReg}
Assume $\epsilon$-regularity (\ref{IneqStandardReg}) holds. There exist constants $C=C(n)$, $\epsilon=\epsilon(n)$ so that if
\be
H&\triangleq&\sup_{B(q,s)\subset{B}(p,r)}\frac{1}{VR(q,s)}\int_{B(q,s)}|\Riem|^{\frac{n}{2}}\label{IneqVolumeEnergy}
\ee
satisfies $H<\epsilon$, then
\be
&&\sup_{B(p,r/2)}|\Riem|\;<\;C\,r^{-2}\,H^\frac2n.\label{IneqHalfRadiusEpsReg}
\ee
\end{proposition}
\underline{\sl Pf} \\
\indent Compare proposition 3.1 of \cite{TV3}, and section 4 of \cite{And3}. We modify the statement of the theorem slightly, proving instead the equivalent statement that $H<\epsilon$ implies $\sup_{B(p,r/4)}|\Riem|\,<\,Cr^{-2}\,H^\frac2n$. By scale-invariance we can assume $r=1$. If all half-radius or smaller subballs, $B(p',s)\subset{B}(p,1)$, $s\le\frac12$ satisfy
\be
&&\sup_{B(p',s/2)}|\Riem|\;<\;Cs^{-2}\left(\epsilon\right)^\frac2n,\label{IneqHalfRadiusEpsRegMod}
\ee
(which is weaker than (\ref{IneqHalfRadiusEpsReg})) then (after perhaps rechoosing $\epsilon$) lemma \ref{LemmaAlmostFlatSobolevBounds} gives that $C_S(B(p,1/2))<C_1(\VR\,B(p,1/2))^{-2/n}$. Therefore
\beg
C_S^{\frac{n}{2}}\int_{B(p,1/2)}|\Riem|^{\frac{n}{2}}&\le&C_1^{\frac{n}{2}}(\VR\,B(p,1/2))^{-1}\int|\Riem|^{\frac{n}{2}}\;\le\;C_1^{\frac{n}{2}}\epsilon_1\;<\;\epsilon_0,
\ees
so that (\ref{IneqStandardReg}) now implies that (\ref{IneqHalfRadiusEpsReg}) holds for $B(p,1/4)$, meaning we get the desired conclusion:
$\sup_{B(p,r/4)}|\Riem|\;<\;Cr^{-2}\,H^\frac2n$

If, however, some half-radius (or smaller) subball does not satisfy (\ref{IneqHalfRadiusEpsRegMod}), then we may replace our original counterexample $B(p,1)$ with a new ball $B(p_1,1/2)$. If again a half-radius (or smaller) subball does not satisfy (\ref{IneqHalfRadiusEpsRegMod}), we get another ball, $B(p_2,1/4)$. Continuing this process, we can find a family of balls $B(p_i,2^{-i})$ so that (\ref{IneqVolumeEnergy}) holds but (\ref{IneqHalfRadiusEpsRegMod}) fails. This process cannot go on indefinitely, because curvature is not infinite on any interior point. Thus, after finitely many steps, we have a ball $B(p_k,2^{-k})$ on which (\ref{IneqHalfRadiusEpsReg}) fails, but so that on any half-radius or smaller ball (\ref{IneqHalfRadiusEpsRegMod}) holds. Rescaling to unit radius, we are in the situation of the previous paragraph, and we get a contradiction. \qed

{\bf Remark}. The main difference, besides the explicitness, between this and proposition 3.1 of \cite{TV3} is that Tian-Viaclovski require comparing the overall energy to each volume ratio, and we require that $\int|\Riem|^{n/2}$ is small compared to the volume ratio on each individual ball.
That said, this theorem is probably most valuable when volume comparison holds, in which case both results are equivalent, and the constant $H$ is controlled by the ratio $\frac{1}{\VR\,B(p,r)}\int_{B(p,r)}|\Riem|^{n/2}$ on just the largest ball.
In the Einstein case one recovers theorem 4.4 of \cite{And3}.

We return to the gradient Ricci soliton case, and establish $\epsilon$-regularity in the form of (\ref{IneqStandardReg}).
In the interest of brevity, and since these arguments appear frequently in the literature, we argue in just the 4-dimensional case (a more complete argument can be found in the references \cite{CS} and \cite{Zhang}).
\begin{proposition}\label{PropSolitonSobolevEpsReg}
There exist constants $C<\infty$, $\epsilon>0$ so that, if $C_S$ is the Sobolev constant for the ball $B(r)$ and
\beg
C_S^{\frac{n}{2}}\cdot\int_{B(r)}|\Riem|^{\frac{n}{2}}\;\le\;\epsilon,
\ees
then
\beg
\sup_{B(r/2)}|\Riem|&\le&C\cdot\left(r^{-2}\,+\,\sup_{B(r)}|\nabla{f}|^2\right)\left(C_S^{\frac{n}{2}}\int_{B(r)}|\Riem|^{\frac{n}{2}}\right)^\frac2n.
\ees
\end{proposition}
\underline{\sl Pf}\\
\indent We scaling the metric so $r=1$; again also we work only with $n=4$.
With $C_S$ denoting the local Sobolev constant and $p\ge2$, the Sobolev inequality gives
\beg
\left(\int\phi^{4}|\Riem|^{2p}\right)^{\frac12}&\le&2C_S\int|\nabla\phi|^p|\Riem|^p\,+\,2C_Sp^2\int\phi^2|\Riem|^{p-2}|\nabla\Riem|^2
\ees
One uses integration-by-parts, a commutator formula, and H\"older's inequality to obtain
\beg
C\int\phi^2|\Riem|^{p-2}|\nabla\Riem|^2&\le&\frac1p\int\phi^2|\Riem|^{p+1}\,+\,\frac{1}{p^2}\int|\nabla\phi|^2|\Riem|^p\,+\,\int\phi^2|\Riem|^{p-2}|\delta\Riem|^2,
\ees
where $C$ is some universal constant. Putting this back in and using proposition \ref{PropEllipticResults},
\beg
C\left(\int\phi^{4}|\Riem|^{2p}\right)^\frac12&\le&C_S\int|\nabla\phi|^2|\Riem|^p\,+\,C_Sp\int\phi^2|\Riem|^{p+1}\,+\,C_Sp^2\int\phi^2|\Riem|^{p-2}|\delta\Riem|^2\\
C\left(\int\phi^{4}|\Riem|^{2p}\right)^\frac12&\le&C_S\int|\nabla\phi|^2|\Riem|^p\,+\,C_Sp\int\phi^2|\Riem|^{p+1}\,+\,C_Sp^2\int\phi^2|\Riem|^p|\nabla{f}|^2
\ees
If $C_S^2\int_{\supp\phi}|\Riem|^2\le\epsilon$, for an appropriate universal $\epsilon$, H\"older's inequality gives
\beg
C\left(\int\phi^{4}|\Riem|^{2p}\right)^\frac12&\le&C_S\,p^2\sup_{\supp\phi}\left(|\nabla\phi|^2\,+\,|\nabla{f}|^2\right)\;\int_{\supp\phi}|\Riem|^p
\ees
It is possible to iterate this estimate. Choosing the $i^{th}$ test function $\phi_i\ge0$ so that $\supp\phi_i\subset{B}(2^{-1}+2^{i-1})$, $\phi_i\equiv1$ on $B(2^{-1}+2^{-i})$, and $\sup|\nabla\phi_i|\le2^{i+2}$, and choosing $p_i=2^{i}$, we get for $i\ge1$
\beg
C\left(\int_{B(2^{-1}+2^{-i-1})}|\Riem|^{2^{i+1}}\right)^{2^{-i-1}}&\le&8^{i2^{-i}}\left(C\cdot{C}_S\cdot\left(1\,+\,\max|\nabla{f}|^2\right)\right)^{2^{-i}}\left(\int_{B(2^{-1}+2^{-i})}|\Riem|^{2^i}\right)^{2^{-i}}.
\ees
So we obtain the standard $\epsilon$-regularity
\beg
\sup_{B(1/2)}|\Riem|&\le&C\cdot\left(1\,+\,\max_{B(1)}|\nabla{f}|^2\right)\left(C_S^2\int_{B(1)}|\Riem|^2\right)^\frac12
\ees
\qed

\begin{theorem}\label{ThmSolitonVolEpsReg1}
Let $(\SSS^n,g,\nabla{f})$ be a gradient Ricci soliton.
There exists a constant $\epsilon(n)>0$ and a constant $C=C(n)<\infty$ so that
\beg
H\triangleq\sup_{B(q,s)\subset{B}(p,r)}\frac{1}{VR\,B(q,s)}\int_{B(q,s)}|\Riem|^2\;\le\;\epsilon
\ees
implies
\beg
\sup_{B(r/2)}|\Riem|&\le&C\cdot\left(r^{-2}\,+\,\sup_{B(p,r)}|\nabla{f}|^2\right)\cdot{H}^{\frac{n}{2}}.
\ees
\end{theorem}
\underline{\sl Pf}\\
\indent Propositions \ref{PropSolitonSobolevEpsReg} and \ref{PropVolEpsReg}. \qed

At finite distances from a basepoint $p_m\in\min_{\SSS}f$, volume ratios are bounded.
Thus we also have

\begin{theorem}\label{ThmSolitonVolEpsReg2}
Let $(\SSS^n,g,\nabla{f})$ be a gradient Ricci soliton.
There exists a constant $\epsilon>0$ and a constant $C<\infty$, both depending on $n$, $\lambda$, $C_1$, and $\dist(p_m,p)$ so that
\beg
\int_{B(r)}|\Riem|^2\;\le\;\epsilon
\ees
implies
\beg
\sup_{B(r/2)}|\Riem|&\le&C\cdot{r}^{-2}\cdot\left(\int_{B(p,r)}|\Riem|^2\right)^\frac12.
\ees
\end{theorem}
\qed

\section{Convergence}\label{SectionConvergence}

This section will be devoted to the proof of our main theorem.
\begin{theorem}
Assume $(\SSS_i,g_i,\nabla{f}_i)\in\mathfrak{S}^n$ is a sequence of n-dimensional unit-volume gradient Ricci solitons. Then a subsequence converges in the pointed Gromov-Hausdorff topology to a unit-volume orbifold gradient Ricci soliton with locally finitely many singular points (in 4 dimensions, finitely many). The singular points are all orbifold points, and away from the singularities the convergence is in the $C^\infty$ topology.
\end{theorem}
\underline{\sl Pf}\\
\indent Let $p_i\in\SSS_i$ be a natural basepoint, choose a large radius $S$, and select an $\epsilon$ so that theorem \ref{ThmSolitonVolEpsReg2} works in the entirety of $B(p_i,S)$.
For any given $r>0$, let $\SSS_{i,r}\subset{B}(p_i,S)$ be the points $q$ with
\be
\int_{B(q,r)}|\Riem|^2\;\le\;\epsilon.\label{IneqNecessaryForRegularity}
\ee
Inside of $B(p_i,S)$ we can make the usual covering argument to bound the number of possible singularities.
Let $G_{i,r,S}$ be the set of points inside $B(p_i,S)$ on which (\ref{IneqNecessaryForRegularity}) fails.
One can cover this set by balls $\{B(q_{ij},2r)\}_{j}$ so that the balls $B(q_{ij},r)$ are disjoint, and $q_{ij}\in{G}_{i,r}$.
The number of balls in such a covering is uniformly finite, since each $B(q_{ij},r)$ has a definite amount of energy.

Put $\SSS_{i,r,S}=B(p_i,S)-\bigcup_j{B}(q_{ij},r)$.
When $r$ is very small, soliton volume comparison implies that most of the volume of $\SSS_i$ is present in $\SSS_{i,r,S}$ (the volume isn't all absorbed into the singularities).
Since sectional curvatures and volume ratios are bounded on $\SSS_{i,r,S}$, a subsequence converges in the $C^{1,\alpha}$ topology, as $i\rightarrow\infty$, to some limiting manifold-with-boundary $\SSS_{\infty,r,S}$.
Sending $r\searrow0$ and passing to further subsequences, we obtain a pointed limiting manifold-with-boundary, whose completion we will denote $(\SSS_{\infty,S},\,p_\infty)$, where $p_\infty$ is the limit of the basepoints $p_i$.
Similarly, letting $S\nearrow\infty$, and passing to further subsequences, we get a limiting manifold with point-type singularities, but otherwise no boundary, the completion of which we denote $(\SSS_\infty,\,p_\infty)$.
Equation (\ref{IneqSuperLevelSetControl}) and the fact that no volume is captured in the bubbles imply that $\Vol\SSS=1$.

Next we consider the regularity of the convergence of the functions $f_i$, and simultaneously of the metric.
We have that $f_i$ is uniformly bounded in balls of uniform size around $p_i$, and therefore scalar curvature $R$ is bounded in these balls.
Thus $\triangle{f}_i=-(R_i+n\lambda_i)$ implies $C^{1,\alpha}$ bounds on the $f_i$, wherever the harmonic radius is controlled.
But then $R_i=C_1-|\nabla{f_i}|^2-2\lambda{f_i}$ gives uniform $C^{0,\alpha}$ bounds on $R_i$, therefore $C^{2,\alpha}$ bounds on $f_i$.
One can continue to bootstrap, but since the ``$\triangle$'' has $g_{ij}$-coefficients, if $g\in{C}^{k,\alpha}$ this bootstrapping terminates with $f\in{C}^{k+2,\alpha}$ and $R\in{C}^{k+1,\alpha}$.
This regularity is uniform at fixed distances from the basepoints $p_i$, and away from any singularities.
Since $C^{k,\alpha}$ bounds on $g$ imply $C^{k+2,\alpha}$ bounds on $f_i$, the equation $\Ric_i=-\nabla^2f_i-\lambda_i{g_i}$ implies that $\Ric$ has $g$ have the same regularity.
Thus wherever it is possible to control the harmonic radius we can use
\beg
\triangle(g_{ij})\;=\;{\Ric}_{ij}\,+\,Q(g,\partial{g})
\ees
to conclude that $g$ has $C^\infty$ bounds.
But the limit occurs locally in the $C^{1,\alpha}$-topology, so construction of harmonic coordinates is possible on the limit manifold, away from the singularities.
The harmonic radius is lower semicontinuous, and so is controlled in the converging sequence as well.
Thus the convergence of $g$ is locally $C^\infty$.

Without going into details, we note that it is possible to bound local Sobolev constants finite distances from the basepoint.
This is due to the noncollapsing of balls, and the availability of relative volume comparison, with which it is possible to carry out the arguments from \cite{Croke} and obtain Sobolev constant bounds.
Then citing the analysis of \cite{CS} \cite{Zhang}, we have that the singularities on $\SSS_\infty$ are of $C^\infty$ orbifold type.

Finally, in dimension 4 we rule out the possibility the limiting orbifold might have infinitely many singularities.
If $q_i\in\SSS_i$ has $\dist(p_i,q_i)$ uniformly bounded and $|\Riem|(q_i)$ unbounded, scale so $|\Riem|(q_i)=1$.
We can assume that, in the rescaled metric, $|\Riem|$ is uniformly bounded on some large ball around $q_i$ (or else we could just rechoose $q_i$; we are capturing a ``deepest bubble'').
We get $C^\infty$ convergence in the pointed Gromov-Hausdorff topology to a complete soliton $X$.
Since $R_i$ and $f_i$ are uniformly bounded a finite distance from $p_i$, the limiting object $X$ has $R=0$, $\lambda=0$, and $C_1=0$ ($\lambda$ and $C_1$ scale like curvature).
Thus the auxiliary equation gives that $|\nabla{f}|^2=0$ so $f$ is constant and $\Ric=0$.
Thus $X$ is an ALE Ricci-flat manifold.
By the splitting theorem, $X$ must be one-ended, or else it would split geometrically as a product of a (Ricci-flat) compact manifold and some $\RR^k$, and therefore be collapsed or be flat $\RR^n$ itself, both of which are impossible.
The Euler number on $X$ is
\beg
\chi(X)&=&\frac{1}{|\Gamma|}\,+\,\frac{1}{8\pi^2}\int|W|^2,
\ees
and is a positive integer.
Here $\Gamma$ is the quotient group at infinity.
If $|\Gamma|=1$, so that the asymptotic volume ratio is exactly Euclidean, the equality case of the Bishop comparison theorem gives that $|\Riem|=0$ everywhere, an impossibility. Thus $\frac{1}{|\Gamma|}\le1/2$, so that $\int|\Riem|^2=\int|W|^2\ge4\pi^2$.

Thus each ``deepest bubble'' must have a uniformly finite amount of energy, so there can only be finitely many of them.  Hence there only finitely many bubbles.
\qed

\noindent{\sl proof of the last statement of corollary \ref{Cor4DimLimits}}\\
\indent Let $M_\infty$ be a $4$-dimensional limiting soliton with basepoint $p_m$.
Let $r_{|R|}(p)=\sup\{r>0\,|\,\sup_{B(p,r)}|\Riem|<r^{-2}\}$ indicate the ``curvature scale'' at $p$.
Define a function $\Sigma(p)=|r_{|R|}(p)|^2(p)\cdot\Vol{B}(p,r_{|R|}(p))$.
We will show that $\Sigma(p)\rightarrow0$ as $\dist(p,p_m)\rightarrow\infty$.
Assuming this is not the case, let $p_i$ be a sequence of points with $\Sigma(p_i)>\epsilon$.

Choose a subsequence $p_i$ so that $\lim_i\Sigma(p_i)=\epsilon$ exists and is nonzero.
We can reselect the points $p_i$ if necessary so that $\Sigma$ `almost' takes on its locally smallest value at $p_i$: if there is some $q_i\in{B}(p_i,2r_{|R|}(p_i))$ with $\Sigma(q_i)<2\Sigma(p_i)$, then reselect setting $p_i=q_i$.
This must terminate after finitely many steps, so that $\Sigma|_{B(p_i,2r_{|R|}(p_i))}\ge2\Sigma(p_i)$.

Scale so $r_{|R|}(p_i)=1$. Then $|\Riem|\le4$ on $B(p,2)$, and we have $\Vol{B}(p_i,1)>\epsilon$, and by Bishop-Gromov volume comparison all subballs of $B(p_i,2)$ have volume ratios controlled by $\frac{1}{100}\epsilon$, say.
But $\int_{B(p_i,1)}|\Riem|^2\searrow0$, so by proposition 2.5 of \cite{And5}, $|\Riem|(p_i)$ must converge to zero on $B(p_i,3/2)$.
But $\sup_{B(p_i,1)}|\Riem|=1$, a contradiction.
\qed

This argument fails in dimension $>4$ because we have not ruled out the possibility of infinitely many singularities, which may accumulate towards infinity.

\end{document}